\documentclass[12pt]{amsart}
\def\ds{\displaystyle}

\newtheorem{theorem}{Theorem}
\title[]
{The Bergman kernel of the symmetrized polydisc in higher dimensions
has zeros}
\author{Nikolai Nikolov and  W\l odzimierz Zwonek}

\address
{Institute of Mathematics and Informatics, Bulgarian Academy of
Sciences, 1113 Sofia, Bulgaria}\email{nik@math.bas.bg}

\address{Instytut Matematyki, Uniwersytet Jagiello\'nski, Reymonta 4,
30-059 Krak\'ow, Poland}\email{Wlodzimierz.Zwonek@im.uj.edu.pl}

\subjclass[2000]{32A25}

\keywords{Bergman kernel, Lu Qi-Keng domain, symmetrized polydisc}

\begin{document}

\begin{thanks}
{This note was prepared during the stay of the first named author
at the Erwin Schr\"odinger Institute, Vienna in September, 2005
and at the Jagiellonian University, Krak\'ow in November, 2005. He
thanks both institutions for their hospitality. The work is a part
of the Research Grant No.1 PO3A 005 28 which is supported by
public means in the programme promoting science in Poland in the
years 2005-2008.}
\end{thanks}

\begin{abstract} We prove that the Bergman kernel of the symmetrized polydisc
in dimension greater than two has zeros.
\end{abstract}

\maketitle

Let $\Bbb D$ be the unit disc in $\Bbb C.$ Let
$\lambda=(\lambda_1,\dots,\lambda_n)\in\Bbb C^n$ and
$\pi_n=(\pi_{n,1},\dots,\pi_{n,n}):\Bbb C^n\to\Bbb C^n$ be defined
as follows:
$$\pi_{n,k}(\lambda)=\sum_{1\le j_1<\dots<j_k\le
n}\lambda_{j_1}\dots,\lambda_{j_k},\ \ 1\le k\le n.$$ The set
$\Bbb G_n=\pi_n(\Bbb D^n)$ is called the symmetrized polydisc. The
symmetrized bidisc $\Bbb G_2$ is the first example of a bounded
pseudoconvex (even hyperconvex) domain that cannot be exhausted by
domains biholomorphic to convex domains and on which the
Carath\'eodory and Kobayashi distances coincide (see
\cite{Cos2003}, \cite{Cos2004} and \cite{Agl-You2004}, see also
\cite{Jar-Pfl}). Note that $\Bbb G_n,$ $n\ge 3,$ cannot be
exhausted by domains biholomorphic to convex domains, too (see
\cite{Nik}). It is, however, not known whether the Carath\'eodory
and Kobayashi distances coincide on $\Bbb G_n$, $n\ge 3$.

In \cite{Edi-Zwo}, the following explicit formula for the Bergman
kernel $K_{\Bbb G_n}$ of $\Bbb G_n$ has been found:
$$K_{\Bbb G_n}(\pi_n(\lambda),\pi_n(\mu))=
\frac{\det[(1-\lambda_j\overline\mu_k)^{-2}]_{1\le j,k\le n}}
{\pi^n\prod_{1\le j<k\le
n}[(\lambda_j-\lambda_k)(\overline\mu_j-\overline\mu_k)]},\ \
\lambda,\mu\in\Bbb D^n.\leqno{(1)}$$ Observe that although the
right-hand side of (1) is not formally defined on the whole $\Bbb
G_n\times\Bbb G_n,$ it extends smoothly on this set. The formula
(1) easily implies that $\Bbb G_2$ is a Lu Qi-Keng domain (see
\cite{Edi-Zwo}), i.e. $K_{\Bbb G_2}$ has no zeros on $\Bbb
G_2\times\Bbb G_2$ -- for the comprehensive information on the Lu
Qi-Keng problem see e.g. \cite{Boas2000}. Then the following
natural question has been posed in \cite{Edi-Zwo} (see also
\cite{Jar-Pfl}): Does $K_{\Bbb G_n}$ have zeros for $n\ge 3$? The
aim of this note is to give a positive answer to the above
question thus providing an example of a proper image of the
polydisc $\Bbb D^n,$ $n\ge 3,$ which is not a Lu Qi-Keng domain.

\begin{theorem}
$K_{\Bbb G_n}$ has zeros for any $n\ge 3.$
\end{theorem}

\begin{proof}

We shall proceed by induction on $n\ge 3$ showing that:

$(\ast)$ there are points $\lambda,\mu\in\Bbb D^n,$ both with
pairwise different coordinates, such that
$$\Delta_n(\lambda,\mu):=\det[(1-\lambda_j\overline\mu_k)^{-2}]
_{1\le j,k\le n}=0$$ and
$f_n:=\Delta_n(\cdot,\lambda_2,\dots,\lambda_n,\mu_1,\dots,\mu_n)\not\equiv
0.$

{\it The case $n=3.$} We have the following formula (see Appendix
A):
$$K_{\Bbb G_3}(\pi_3(\lambda_1,\lambda_2,\lambda_3),\pi_3(\mu_1,\mu_2,0))=
\frac{a(\nu)z^2-b(\nu)z+2c(\nu)}{\pi^3\prod_{1\le j\le 3,1\le k\le
2}(1-\lambda_j\overline\mu_k)^2},\leqno{(2)}$$ where $\ds
z=\frac{\overline\mu_2}{\overline\mu_1}$ $(\mu_1\neq 0),$
$\nu_j=\lambda_j\overline\mu_1,$ $j=1,2,3,$ and
$$a(\nu)=\pi_{3,2}(\nu)(2-\pi_{3,1}(\nu))+\pi_{3,3}(\nu)(2\pi_{3,1}(\nu)-3),$$
$$b(\nu)=(\pi_{3,1}(\nu)-2)(\pi_{3,2}(\nu)-2\pi_{3,1}(\nu)+3)+
3(\pi_{3,3}(\nu)-\pi_{3,1}(\nu)+2),$$
$$c(\nu)=\pi_{3,2}(\nu)-2\pi_{3,1}(\nu)+3.$$
For the fixed point $\nu_0=(e^{i\pi/6},e^{i\pi/3},e^{-i\pi/6})
\footnote{We thank Pencho Marinov for his computer programme which
shortened the calculations needed for finding good candidates for
zeros of $K_{\Bbb G_3}$.}$ the number
$$z_0=e^{-i\pi/4}\frac{6-3\sqrt3-\sqrt{40\sqrt3-69}}{\sqrt2(3\sqrt3-5)}$$
satisfies the equality $a(\nu_0)z_0^2-b(\nu_0)z_0+2c(\nu_0)=0$
(see Appendix B). Since $z_0\in\Bbb D,$ it follows that for any
$\nu\in\Bbb D^3,$ close to $\nu_0,$ there is a $z\in\Bbb D,$ close
to $z_0,$ such that $a(\nu)z^2-b(\nu)z+2c(\nu)=0.$ Then choosing
$\mu_1\in\Bbb D$ with $|\mu_1|>|\nu_1|,|\nu_2|,|\nu_3|$ we get
points $\lambda,\mu\in\Bbb D^3,$ both with pairwise different
coordinates such that $\Delta_3(\lambda,\mu)=0.$

To see that $f_3\not\equiv 0$ assume the contrary. Then
$f_3(0)=f_3'(0)=f_3''(0)=0,$ i.e.
$$\det\left[\begin{matrix}
\overline\mu_1^j&\overline\mu_2^j&\overline\mu_3^j\\
(1-\lambda_2\overline\mu_1)^{-2}&(1-\lambda_2\overline\mu_2)^{-2}&(1-\lambda_2\overline\mu_2)^{-2}\\
(1-\lambda_3\overline\mu_1)^{-2}&(1-\lambda_3\overline\mu_2)^{-2}&(1-\lambda_3\overline\mu_3)^{-2}\\
\end{matrix}\right]=0$$
for $j=0,1,2.$ Since $\mu_1,\mu_2,\mu_3$ are pairwise different,
the vectors $(1,1,1),$ $(\mu_1,\mu_2,\mu_3)$ and
$(\mu_1^2,\mu_2^2,\mu_3^2)$ are linearly independent. It follows
that the vectors in the second and the third lines of the above
determinant are linearly dependent. In particular, $K_{\Bbb
G_2}(\pi_2(\lambda_2,\lambda_3),\pi_2(\mu_2,\mu_3))=0,$ a
contradiction.

{\it The induction step.} Assume that $(\ast)$ holds for some $n\ge
3.$ We shall choose numbers $\tilde\lambda_1$ and
$\tilde\lambda_{n+1}$ in $\Bbb D,$ close to $\lambda_1$ and 1,
respectively (which will provide pairwise different coordinates of
the new points), such that
$$g_{n+1}(\tilde\lambda_1,\tilde\lambda_{n+1}):=
\Delta_{n+1}(\tilde\lambda_1,\lambda_2,\dots,\lambda_n,\tilde\lambda_{n+1},
\mu_1,\dots,\mu_n,\tilde\lambda_{n+1})=0$$ and
$g_{n+1}(\cdot,\lambda_{n+1})\not\equiv 0.$ Note that
$$g_{n+1}(\tilde\lambda_1,\tilde\lambda_{n+1})=\frac{f_n(\tilde\lambda_1)}
{(1-|\tilde\lambda_{n+1}|^2)^2}+h_n(\tilde\lambda_1,\tilde\lambda_{n+1}),$$
where $h_n$ is a continuous function on $\Bbb
D\times\overline{\Bbb D}.$ Since $f_n\not\equiv 0$ is a
holomorphic function, for any small $r>0$ the number $\lambda_1$
is the only zero of $f_n$ in the closed disc $D\subset\Bbb D$ with
center at $\lambda_1$ and radius $r.$ Then
$m:=\ds\frac{\min_{\partial D}|f_n|}{\max_{\partial D
\times\overline{\Bbb D}}|h_n|}>0.$ Hence
$|f_n|>(1-|\tilde\lambda_{n+1}|^2)^2|h_n(\cdot,\tilde\lambda_{n+1})|$
on $\partial D$ if $1-|\tilde\lambda_{n+1}|^2<\sqrt m.$  Fix such
a $\tilde\lambda_{n+1}$. Since $h_n(\cdot,\tilde\lambda_{n+1})$ is
a holomorphic function on $\Bbb D$, the Rouch\'e theorem implies
that $g_{n+1}(\cdot,\tilde\lambda_{n+1})$ has finitely many zeros
in $D$ (in particular,
$g_{n+1}(\cdot,\tilde\lambda_{n+1})\not\equiv 0$), which completes
the proof.

\end{proof}

\noindent{\it Remark.} The above proof shows that if $n\ge 4,$
then there are points  $(\lambda,\nu),$ close to the diagonal of
$\Bbb D^n\times\Bbb D^n$ in the following sense:
$\lambda_j=\mu_j>0$ for $j=4,\dots,n$ such that $K_{\Bbb
G_n}(\pi_n(\lambda),\pi_n(\mu))=0.$ On the other hand, it can be
shown that $K_{\Bbb G_3}(\pi_3(\lambda),\pi_3(\mu))\neq0$ if
$\lambda_3=\mu_3.$

\

\noindent{\bf Appendix A.} By (1), one has that
$$\pi^3(\lambda_1-\lambda_2)(\lambda_1-\lambda_3)(\lambda_2-\lambda_3)\overline\mu_1
\overline\mu_2(\overline\mu_1-\overline\mu_2)K_{\Bbb
G_3}(\pi_3(\lambda_1,\lambda_1,\lambda_3),\pi_3(\mu_1,\mu_2,0))\leqno{(3)}$$
$$=\det\left[\begin{matrix}
(1-\nu_1)^{-2}&(1-z\nu_1)^{-2}&1\\
(1-\nu_2)^{-2}&(1-z\nu_2)^{-2}&1\\
(1-\nu_3)^{-2}&(1-z\nu_3)^{-2}&1
\end{matrix}\right]$$
$$=\det\left[\begin{matrix}
(1-\nu_1)^{-2}-(1-\nu_3)^{-2}&(1-z\nu_1)^{-2}-(1-z\nu_3)^{-2}\\
(1-\nu_2)^{-2}-(1-\nu_3)^{-2}&(1-z\nu_2)^{-2}-(1-z\nu_3)^{-2}
\end{matrix}\right]$$
$$=\frac{(\nu_1-\nu_3)(\nu_2-\nu_3)z}{(1-\nu_3)^2(1-z\nu_3)^2}
\det\left[\begin{matrix}\ds
\frac{\nu_1+\nu_3-2}{(1-\nu_1)^2}&\ds\frac{z\nu_1+z\nu_3-2}{(1-z\nu_1)^2}\\
\ds\frac{\nu_2+\nu_3-2}{(1-\nu_2)^2}&\ds\frac{z\nu_2+z\nu_3-2}{(1-z\nu_2)^2}
\end{matrix}\right]$$
$$=\frac{(\nu_1-\nu_3)(\nu_2-\nu_3)z}{\prod_{1\le j\le 3,1\le
k\le 2}(1-\lambda_j\overline\mu_k)^2}
\Bigl((\nu_1+\nu_3-2)(z\nu_2+z\nu_3-2)(1-z\nu_1)^2(1-\nu_2)^2$$
$$-(\nu_2+\nu_3-2)(z\nu_1+z\nu_3-2)(1-\nu_1)^2(1-z\nu_2)^2\Bigr)\leqno{(4)}$$
$$=\frac{(\nu_1-\nu_3)(\nu_2-\nu_3)z(z-1)(A(\nu)z^2-B(\nu)z+2C(\nu))}{\prod_{1\le j\le 3,1\le
k\le 2}(1-\lambda_j\overline\mu_k)^2}.\leqno{(5)}$$ To find
$A(\nu)$, $B(\nu)$ and $C(\nu)$, we shall use that the
coefficients of $z^3,$ $z^0$ and $z$ in the large brackets in (4)
are equal to
$$A(\nu)=(\nu_1+\nu_3-2)(\nu_2+\nu_3)\nu_1^2(1-\nu_2)^2-
(\nu_2+\nu_3-2)(\nu_1+\nu_3)\nu_2^2(1-\nu_1)^2,$$
$$-2C(\nu)=2(\nu_2+\nu_3-2)(1-\nu_1)^2-2(\nu_1+\nu_3-2)(1-\nu_2)^2\hbox{ and}$$
$$B(\nu)+2C(\nu)=(\nu_1+\nu_3-2)(\nu_2+\nu_3+4\nu_1)(1-\nu_2)^2$$
$$-(\nu_2+\nu_3-2)(\nu_1+\nu_3+4\nu_2)(1-\nu_1)^2,$$ respectively.
Calculations lead to the formulas
$$A(\nu)=(\nu_2-\nu_1)(\pi_{3,2}(\nu)(2-\pi_{3,1}(\nu))+\pi_{3,3}(\nu)
(2\pi_{3,1}(\nu)-3)),$$
$$C(\nu)=(\nu_2-\nu_1)(\pi_{3,2}(\nu)-2\pi_{3,1}(\nu)+3),$$
$$B(\nu)=(\nu_2-\nu_1)((\pi_{3,1}(\nu)-2)(\pi_{3,2}(\nu)-2\pi_{3,1}(\nu)+3)
+3(\pi_{3,3}(\nu)-\pi_{3,1}(\nu)+2)).$$ To get (2), it remains to
substitute these formulas in (5) and then to compare (5) and (3).

\

\noindent{\bf Appendix B.} Since
$$\pi_{3,1}(\nu_0)=\frac{1+2\sqrt3+i\sqrt3}{2},\
\pi_{3,2}(\nu_0)=\frac{2+\sqrt3+i3}{2},\
\pi_{3,3}(\nu_0)=e^{i\pi/3},$$ the formulas for $a(\nu),$ $b(\nu)$
and $c(\nu)$ lead to
$$a(\nu_0)=(3\sqrt3-5)e^{i\pi/3},\ b(\nu_0)=
(6\sqrt2-3\sqrt6)e^{i\pi/12},\ c(\nu_0)=(2\sqrt3-3)e^{-i\pi/6}.$$
Then for $z=e^{-i\pi/4}x$ one has
$e^{i\pi/6}(a(\nu_0)z^2-b(\nu_0)z+2c(\nu_0))$ $$=
(3\sqrt3-5)x^2+(3\sqrt6-6\sqrt2)x+4\sqrt3-6=:p(x).$$ The zeros of
the polynomial $p$ are equal to
$\ds\frac{6-3\sqrt3\pm\sqrt{40\sqrt3-69}}{\sqrt2(3\sqrt3-5)}.$
Note that the smaller one lies in $(0,1),$ since $p(0)>0>p(1).$

\end{document}